\documentclass[reqno,12pt]{amsart}
\usepackage{epsf}
\usepackage{todonotes}
\usepackage{verbatim}
\usepackage{amsmath}
\usepackage{amsopn}
\usepackage{amssymb}
\usepackage{latexsym}
\usepackage{amsfonts}
\usepackage{amsthm}
\usepackage{tikz}

\newcommand{\R}{{\mathbb{R}}}

\newtheorem{thm}{Theorem}[section]

\theoremstyle{definition}

\newtheorem{exa}[thm]{Example}

\newtheorem{defn}[thm]{Definition}

\begin{document}

\title[Signed partitions - A balls into urns approach]{Signed partitions -\\ A balls into urns approach}

\thanks{This research was supported by a grant from the Ministry of Science and Technology, Israel, and the France's Centre National pour la Recherche Scientifique (CNRS)}

\author{Eli Bagno and David Garber}

\address{Eli Bagno, Jerusalem College of Technology\\
21 Havaad Haleumi St. Jerusalem, Israel}
\email{bagnoe@g.jct.ac.il}

\address{David Garber \\
Department of  Applied Mathematics, Holon Institute of Technology,
52 Golomb St., PO Box 305, 58102 Holon, Israel,  and (sabbatical:) Einstein Institute of Mathematics, Hebrew University of Jerusalem, Jerusalem, Israel }
\email{garber@hit.ac.il}

\begin{abstract}
Using Reiner's definition of  Stirling numbers of type B of the second kind, we provide a 'balls into urns' approach for proving a generalization of a well-known identity concerning the classical Stirling numbers of the second kind:
$x^n=\sum\limits_{k=0}^n{S(n,k)[x]_k}.$
\end{abstract}


\maketitle

\section{Introduction}
The partitions of the set $[n]=\{1,\dots,n\}$ in $k$ blocks are enumerated by the {\em Stirling numbers of the second kind}, denoted by $S(n,k)$
(see \cite[page 81]{EC1}). These numbers arise in a variety of problems in enumerative combinatorics; they have many combinatorial interpretations, and have been generalized in various contexts and in many ways.

One of the celebrated results concerning Stirling numbers of the second kind is the following:
\label{thm:typeA_falling}
Let $x \in \mathbb{R}$ and let $n \in \mathbb{N}$. Then we have:
\begin{equation}\label{reg}
x^n=\sum\limits_{k=0}^n{S(n,k)[x]_k},
\end{equation}
\noindent
where $[x]_k:=x(x-1) \cdots (x-k+1)$ is the {\it falling factorial of degree $k$} and $[x]_0:=1$.

This identity arises when one expresses the standard basis of the polynomial ring $\R[x]$ as a linear combination of the basis consisting of the falling factorials (see e.g. the survey of Boyadzhiev \cite{Boy}).

\medskip

There are some known proofs for this identity. A combinatorial one, realizing $x^n$ as the number of functions from the set  $\{1,\dots,n \}$ to the set $\{1,\dots,x \}$ (for an integer $x$), is presented by Stanley \cite[Eqn. (1.94d); its proof is in page 83]{EC1}, and we quote it here (where $\#N = n$ and\break $\#X = x$):

\medskip

\hangindent2.5em
\hangafter=0
\noindent
``The left-hand side is the total
number of functions $f : N \to X$. Each such function is surjective onto a unique subset
$Y = f(N)$ of $X$ satisfying $\#Y \leq n$. If $\#Y = k$, then there are $k!S(n, k)$ such functions, and
there are ${x \choose k}$ choices of subsets $Y$ of $X$ with $\#Y = k$. Hence:\\
\centerline{$x^n =\sum\limits_{k=0}^n k!S(n, k) {x \choose k} = \sum\limits_{k=0}^n
S(n, k)[x]_k.$''}

\medskip

There is a nice generalization of Identity (\ref{thm:typeA_falling}), which appears in Remmel and Wachs \cite{RW} and Bala \cite{Bala}.
In order to demonstrate this generalization combinatorially, we use the Stirling numbers of type $B$ of the second kind, denoted by $S_B(n,k)$, which are related to the Coxeter group of type B. The exact definition will be given in Section \ref{pre}. Their generalization is:

\begin{thm}\label{thm_bala}
Let $x \in \mathbb{R}$ and let $n \in \mathbb{N}$. Then we have:
\begin{equation}
x^n=\sum\limits_{k=0}^n{S_B(n,k)[x]^B_k},\label{B}
\end{equation}
where $[x]^B_k:=(x-1)(x-3)\cdots (x-2k+1)$ and $[x]^B_0:=1$.
\end{thm}

Remmel and Wachs \cite{RW} proved this equality by using the combinatorial interpretation of $S_{B}(n,k)$ as counting configurations of $k$-non attacking rooks (specifically, this is $S_{n,k}^{0,2}(1,1)$ in their notation).  Bala \cite{Bala} proved this equality using a generating-functions technique ($S_{(2,0,1)}$ in his notation).

In \cite{BBG}, a geometric way to prove Equation \eqref{B}, interpreting $x^n$ as counting the number of points in a cubical lattice, is presented.
\medskip

The purpose of this note (see  Section \ref{proof}) is a simple combinatorial proof, which interprets both sides of Equation \eqref{B}, using a balls into urns approach.  Note that our proof is actually a generalization for Coxeter groups of type $B$ of the proof  for Equation \eqref{thm:typeA_falling}, that we have quoted above.

\section{Signed partitions}\label{pre}
We define  the objects which the Stirling numbers of type $B$ of the second kind count, introduced by Reiner \cite{R}. Denote $[\pm n]:=\{\pm1,\dots,\pm n\}$ and
$-C:=\{-i \ | \ i\in C\}$ for a set $C \subseteq [\pm n]$.

\begin{defn}
A {\it signed partition} is a set partition of $[\pm n]$ into blocks, which satisfies the following conditions:
\begin{itemize}
\item There exists at most one block satisfying $-C=C$, called the {\em zero-block}. It is a subset of $[\pm n]$ of the form $\{\pm i \mid i \in S\}$ for some $S \subseteq [n]$.
\item If $C$ appears as a block in the partition, then $-C$ also appears in that partition.
\end{itemize}
\end{defn}

We denote by $S_B(n,k)$ the number of  signed partitions of $[\pm n]$ having exactly $k$ pairs of nonzero blocks. These numbers are called {\em Stirling numbers of type $B$ of the second kind}. They form the sequence A039755 in the OEIS \cite{oeis}. Table~\ref{Table1} records these numbers for small values of $n$ and $k$.

\begin{table}[!ht]
\begin{center}
\begin{tabular}{r||r|r|r|r|r|r|r}
  $n/k$ & 0 & 1 & 2 & 3 & 4 & 5 & 6  \\
  \hline\hline
   0    & 1 &   &   &   &   &   &     \\
   1    & 1 & 1 &   &   &   &   &     \\
   2    & 1 & 4 & 1  &   &   &   &     \\
   3    & 1 & 13 & 9  & 1  &   &   &     \\
   4    & 1 & 40 & 58  & 16  & 1  &   &     \\
   5    & 1 & 121 & 330  & 170  & 25  & 1  &     \\
   6    & 1 & 364 & 1771  & 1520  & 395  & 36  & 1    \\
\end{tabular}
\end{center}

\caption{Stirling numbers of type $B$ of the second kind $S_B(n,k)$.}\label{Table1}
\end{table}

\begin{exa} The following partitions
$$P_1=\{ \{ 3, -3\}, \{ -2,1 \},\{2,-1\},\{ -4,5 \},\{4,-5\}\},$$
$$P_2=\{ \{ 3\}, \{ - 3\},\{ -2,1 \},\{ 2,-1 \}, \{ -4,5 \},\{ 4,-5 \} \},$$
are respectively a signed partition of $[\pm 5]$ with a zero block $\{ 3, -3\}$ and a signed partition of $[\pm 5]$ without a zero-block.
\end{exa}

\section{The combinatorial proof}\label{Falling polynomials for Coxeter groups}\label{proof}

In this section, we supply a direct combinatorial proof for Theorem \ref{thm_bala}, where $x^n$ is interpreted as the number of assignments of $n$ balls numbered $1$ to $n$ into $x$ distinguishable urns. As will be explained below, we can assume that $x$ is an integer.


\begin{proof}[Direct combinatorial proof]
Since Equation (\ref{B}) is a polynomial identity, it is sufficient to prove it for odd natural numbers. Let $m \in \mathbb{N}$ be an odd number. We will show:
\begin{equation} \label{eq}
m^n=\sum\limits_{k=0}^n{S_B(n,k)[m]^B_k}.
\end{equation}
The left-hand side of Equation (\ref{eq}) is the number of assignments of $n$ balls numbered $1$ to $n$  into $m$ urns.  In the right-hand side, we associate $$[m]^B_k =(m-1)(m-3)\cdots (m-2k+1)$$ assignments to each one of the $S_B(n,k)$ signed partitions, and then we sum them up to get the total number of assignments, thus proving the identity.

\medskip

Let $\mathcal{B} = \{B_0,B_1,-B_1,\dots,B_k , -B_k\}$
be a signed partition, where $B_0$ is the zero-block (which possibly does not exist).
Note that by the definition of $[x]_k^B$ which appeared in Theorem \ref{thm_bala}, we have $[x]_k^B=0$ for $m<2k$, so we may assume that $k < \frac{m}{2}$.

For our convenience, we impose the following order on the blocks of the signed partition: the blocks $B$ and $-B$ are adjacent and the pairs of blocks of $\mathcal{B}$ are ordered in such a way that pairs of blocks which have smaller minimal positive elements precede (except for the zero-block $B_0$ which is always located as the first block). For each pair of blocks $B$ and $-B$, the internal order between
$B$ and $-B$ is chosen in such a way that the block which contains the minimal positive element of $B\cup -B$ is located first.
For example, $$\{\{5,-5\},\{1,-3\},\{-1,3\},\{2,4\},\{-2,-4\}\}$$
is properly ordered.

\medskip

For convenience, we consider an assignment of $n$ balls into $m$ urns as a function $f: [n]\rightarrow [m]$,
and associate with $\mathcal{B}$ the set of ball assignments according to the following procedure:

\begin{itemize}

\item For any {\it positive} $i\in B_0$, define: $f(i)=1$.
\item Choose a number $p$ out of the $m-1$ remaining numbers\break ($2 \leq p \leq m)$, and send the positive elements of $B_1$ to $p$. The absolute values of the negative elements of $B_1$ (i.e. the positive elements of $-B_1$, if they exist) will be sent to the next number in cyclical order excluding the number $1$ (which might have already been occupied by the positive elements of the zero-block). This can be done in $m-1$ different ways.

\item We pass to the pair of blocks $B_2$ and $-B_2$. Similarly, choose a new number $p'$ out of the $m-3$ remaining numbers (the number $1$ is occupied by the positive elements of the zero-block, and two additional numbers are already occupied by the elements of the pair of blocks $B_1$ and $-B_1$), and send the positive elements of $B_2$ to $p'$. For each negative $i \in B_2$, the absolute value of $i$ will be sent to the next unoccupied number in cyclical order. This may be done in $m-3$ different ways.

\item Proceeding this way, we associate a set of $[m]^B_k$ functions from $[n]$ to $[m]$ to each signed partition having $k$ pairs of nonzero blocks.

\end{itemize}

\medskip

Conversely, we now recover the signed partition from a given function $f:[n] \rightarrow [m]$. Define: $$B_0=\{\pm i \mid f(i)=1\}.$$ Mark the number $1$ as used.  Let $k \in [n]$ be the minimal positive number such that $f(k)\neq 1$. Denote $a:=f(k)$.
Let $b \in [m]-\{1\}$ be the next unused number in cyclical order. Define:
$$B_1=\{i \in [n] \mid f(i)=a\} \cup \{-i \mid f(i)=b\},$$ and add the blocks $B_1$ and $-B_1$ to the signed partition.  Now mark the numbers $a,b$ as used.
Proceeding along these lines, we arrive at the signed partition which induces the function $f$.
\end{proof}

\medskip

\begin{exa}\label{exam3.2}
Let $n=6$ and $m=7$. Consider the signed partition
$$\mathcal{B}=\{\{\pm 1\},\{2,-3,5\},\{-2,3,-5\},\{4,-6\},\{-4,6\}\}$$
of $[\pm 6]$. Every function $f:[6] \rightarrow [7]$ which is induced by $\mathcal{B}$ sends $1$ (which is the content of the zero block) to $1$ (see Figure \ref{fig1}).

\begin{figure}[!ht]
    \centering
\begin{tikzpicture}[cap=round,line width=2pt]

\foreach \i in
    {1, 2, 3, 4, 5, 6,7}
  {
    \draw[-] (2*\i-2,0) -- (1+2*\i-2,0);
    \draw[-] (2*\i-2,0) -- (0+2*\i-2,1.8);
    \draw[-] (2*\i-1,0) -- (2*\i-1,1.8);
    \node[below] at (0.5+2*\i-2,0){\i};
  }

\draw[line width=1pt]  (0.5,0.5) circle (0.3cm);

\node at (0.5,0.5){1};

\draw[line width=1pt]  (6.5,0.5) circle (0.3cm);

\node at (6.5,0.5){2};

\draw[line width=1pt]  (6.5,1.3) circle (0.3cm);

\node at (6.5,1.3){5};

\draw[line width=1pt]  (8.5,0.5) circle (0.3cm);

\node at (8.5,0.5){3};

\draw[line width=1pt]  (12.5,0.5) circle (0.3cm);

\node at (12.5,0.5){4};

\draw[line width=1pt]  (2.5,0.5) circle (0.3cm);

\node at (2.5,0.5){6};

\draw[line width=1pt, dashed] (13.4,0) .. controls
(12.8,-0.9) and (12.3,-0.95)  ..
(12,-1);
\draw[line width=1pt, dashed] (12,-1) -- (3,-1);
\draw[line width=1pt,dashed,->] (3,-1) .. controls (2.7,-0.95)  and (2.2,-0.9).. (1.5,0);

\end{tikzpicture}
\caption{An assignment of $6$ balls into $7$ urns}\label{fig1}
\end{figure}
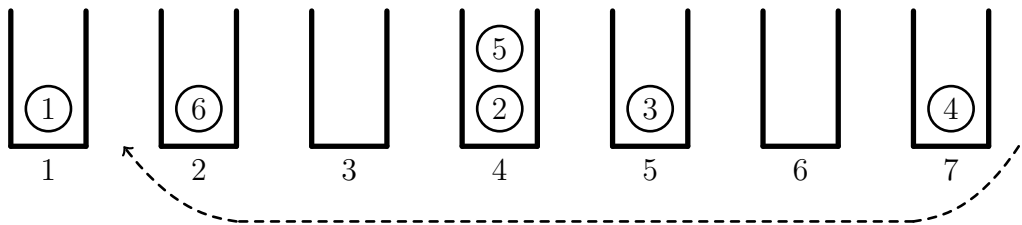

Now we pass to the first block $\{2,-3,5\}$ (together with its negative block $\{-2,3,-5\}$): we have to choose a value for the images of $2$ and $5$ out of $6$ possibilities. Take for example $f(2)=f(5)=4$. Then we have to assign $f(3)=5$, which is the next free value in cyclical order.

The next block is $\{4,-6\}$ (together with its negative block $\{-4,6\}$). For this block we are left with $4$ possibilities for assigning values. Choose for instance $f(4)=7$ and so we must assign $f(6)=2$, which is the next free value in cyclical order. The resulting balls into urns assignment is depicted in Figure \ref{fig1}.

\medskip

Conversely, given the balls into urns assignment obtained above:
$$f(1)=1,\ f(2)=4,\ f(3)=5,\  f(4)=7,\  f(5)=4,f(6)=2.$$
In order to recover the signed partition $\mathcal{B}$ which induced this assignment, we act as follows:

\begin{itemize}
\item Only $1$ is sent to $1$, so we have the zero block $B_0=\{\pm 1\}$.
\item The current minimal unused element is $2$ which is sent by $f$ to $4$, so the positive part of the next block will be  $f^{-1}(\{4\})=\{2,5\}$ and the negative part will be $f^{-1}(\{5\})=\{3\}$ (since $5$ is the next value after $4$ in cyclical order). Hence we get the pair of blocks: $\{2,5,-3\}$ and $\{-2,-5,3\}$.
\item Now, the current minimal unused element is $4$
which is sent by $f$ to $7$, so the positive part of the block will be  $f^{-1}(\{7\})=\{4\}$ and the negative part will be $f^{-1}(\{2\})=\{6\}$ (since $2$ is the next value after $7$ in cyclical order). Hence we get the pair of blocks: $\{4,-6\}$ and $\{6,-4\}$.
\end{itemize}

\medskip

So we get that the signed partition is:
$$\mathcal{B}=\{\{\pm 1\},\{2,5,-3\},\{-2,-5,3\},\{4,-6\},\{-4,6\}\},$$
which indeed was our original signed partition.
\end{exa}

\end{document}